\documentclass[11pt]{article}
\usepackage{color}
\usepackage[colorlinks=true]{hyperref}
\usepackage[margin=1in]{geometry}
\usepackage{amssymb}
\usepackage{amsmath}
\usepackage{amsthm}
\usepackage{mathtools}
\usepackage{epsfig, graphicx}
\newtheorem{theorem}{Theorem}

\newtheorem{remark}{Remark}

\allowdisplaybreaks   % For long formula
\begin{document}
\title{Solving High Dimensional Partial Differential Equations Using Tensor Type Discretization and
Optimization Process\footnote{This work was supported in part by the National Key Research and
Development Program of China (2019YFA0709601), the National Center for Mathematics and Interdisciplinary
Science, CAS.}}
\author{Yangfei Liao\footnote{LSEC, NCMIS, Institute
of Computational Mathematics, Academy of Mathematics and Systems
Science, Chinese Academy of Sciences, Beijing 100190,
China,  and School of Mathematical Sciences, University
of Chinese Academy of Sciences, Beijing 100049, China (yangfeiliao@lsec.cc.ac.cn).} ,\ \
Yifan Wang\footnote{LSEC, NCMIS, Institute
of Computational Mathematics, Academy of Mathematics and Systems
Science, Chinese Academy of Sciences, Beijing 100190,
China,  and School of Mathematical Sciences, University
of Chinese Academy of Sciences, Beijing 100049, China (wangyifan@lsec.cc.ac.cn).} \ \
 and \ \
Hehu Xie\footnote{LSEC, NCMIS, Institute
of Computational Mathematics, Academy of Mathematics and Systems
Science, Chinese Academy of Sciences, Beijing 100190,
China,  and School of Mathematical Sciences, University
of Chinese Academy of Sciences, Beijing 100049, China (hhxie@lsec.cc.ac.cn).}}
%=====================================================================================
\date{}
\maketitle
%=====================================================================================
\begin{abstract}
In this paper,  we propose a tensor type of discretization and optimization process for solving high dimensional
partial differential equations.  First, we design the tensor type of trial function for the high dimensional
partial differential equations. Based on the tensor structure of the trial functions, we can do the direct numerical integration of 
the approximate solution without the help of Monte-Carlo method.  Then combined with the Ritz or Galerkin
method, solving the high dimensional partial differential equation can be transformed to solve a concerned
optimization problem. Some numerical tests are provided to validate the proposed numerical methods.
%-----------------------------------------------------------------------------------------------------------------------------------
\end{abstract}
%=====================================================================================

\section{Introduction}
In this paper, we design a type of numerical method to solve high dimensional partial differential equations 
by using the standard finite element method and the optimization process.
The key point is to use the tensor decomposition structure for representing the high dimensional trial function with some type of one dimensional basis functions. Solving partial differential equations is a basic and important task in many scientific and industrial applications. There also exist many high-dimensional PDEs such as  many-body Schr\"{o}dinger, Boltzmann equations,  Fokker-Planck equations, stochastic PDEs (SPDEs), which are almost impossible solved by the traditional numerical methods.

Recently, the neural networks (NNs) based machine learning method is a popular way to
solve the high-dimensional PDEs  (\cite{BaymaniEffati,EYu,LagarisLikasPapageorgiou,RaissiPerdikarisKarniadakis,SirignanoSpiliopoulos}).
The reason is that NNs can approximate any function if it is given enough parameters.
The fully-connected neural network (FNN) is the most widely used architecture to build
the functions for solving high-dimensional PDEs since its universal approximation property.
So far, there are several types of FNN-based methods such as well-known deep Ritz \cite{EYu},
PINN \cite{RaissiPerdikarisKarniadakis} and  DGM \cite{SirignanoSpiliopoulos} 
% and weak adversarial network \cite{ZangBaoYeZhou},
for solving high-dimensional PDEs by designing different types of loss functions.
Among these methods, the loss functions always include computing high-dimensional integrations
for the functions defined by FNN.
Direct numerical integration for the high-dimensional functions always meets the ``curse of dimensionality'' (CoD).
Then Monte-Carlo method is always adopted to compute these high-dimensional integration with
some types of sampling methods \cite{EYu}.

This paper proposes  a new type of tensor type of discretization method to solve high dimensional
partial differential equations without the help of  Monte-Carlo process.
The CANDECOMP/PARAFAC (CP) tensor decomposition \cite{BeylkinMohlenkampInt,HongKoldaDuersch,KoldaBader} 
can build the low-rank approximation to high
dimensional functions.  The CP tensor decomposes can be considered as the higher-order extensions of
the singular value decomposition (SVD) for the matrices. This means the SVD provides an idea to
decompose the high-dimensional Hilbert space into the tensor product of several Hilbert spaces \cite{Ryan}.
The tensor product decomposition has been used to establish low-rank approximations of operators and functions \cite{BeylkinMohlenkamp,HackbuschKhoromskij,ReynoldsDoostanBeylkin}.

%It is worth to mentioning that there is no known general results to give the relationship between the rank and error bounds.
%For more details, please refer to \cite{HongKoldaDuersch,KoldaBader} and numerical investigations \cite{BeylkinMohlenkamp}.

%In this paper, we use the CP decomposition way to build the trial functions directly.  

The aim of this paper is to propose a type of tensor type discretization (TTD) method to solve
high-dimensional PDEs, which is based on the direct CP decomposition way to represent the trial 
functions for the high dimensional partial differential equations. 
With the help of some type of one dimensional basis functions, 
in the TTD, the trial function is constructed with the CP tensor decomposition way.
Based on the tensor type of structure, we can design a efficient numerical integration
for high dimensional trial functions.
We will show, the computational work for the integration of these functions is only linear scale of
the dimension. Based on this type of efficient numerical integrations, solving high dimensional partial differential
equations can be transformed to the corresponding optimization problem. 
It is worth to mentioning that there is no tensor decomposition in TTD. We only use the tensor decomposition 
to build (or represent) the trial functions and then solve the deduced optimization problems with some type of 
iterative method. Finally, the TTD  provides a way to overcome CoD in some sense
for solving high-dimensional PDEs, which is the main motivation and contribution of this paper.

An outline of the paper goes as follows. In Section \ref{Section_TTD}, we introduce the tensor type of method to build the high dimensional trial functions.
Section \ref{Section_BVP} is devoted to proposing the TTD method for solving the high-dimensional 
boundary value problem with the efficient numerical integration for the trial functions.
Some numerical examples are provided in Section \ref{Section_Numerical}
to show the validity and efficiency of the proposed numerical methods in this paper.
Some concluding remarks are given in the last section.

%-----------------------------------------------------------------------------------------------------
\section{Tensor type of trial function}\label{Section_TTD}
In the section, we introduce the tensor type of way to build the trial function to approximate the
solution of the high dimensional partial differential equations. The aim here is to build some type of
low rank approximations for the high dimensional solutions.

In order to describe the way to build the high dimensional trail function to approximate the solution of
partial differential equations,
for $i=1, \cdots, d$, we use $\{\phi_{i,1}(x_i), \cdots, \phi_{i,N_i}(x_i)\}$ to denote the basis of
the finite dimensional space $\mathcal V_i$  which is a subspace of $H_0^1(\Omega_i)$ on the one dimensional domain.
Based on these basis, we can construct the discrete space $\mathcal V \subset H_0^1(\Omega)$
by using the tensor product way, i.e.,
$$\mathcal V:= \mathcal V_1 \otimes \mathcal V_2\otimes\cdots\otimes\mathcal V_d \subset H_0^1(\Omega_1\times\cdots\times\Omega_d) = H_0^1(\Omega).$$
With the help of tensor product way, we can also construct the basis as follows
\begin{eqnarray}
\mathcal V={\rm span}\left\{ \prod_{i=1}^d\phi_{i,j_i}(x_i),\ \ 1\leq j_i \leq N_i\right\}.
\end{eqnarray}
It is obvious that ${\rm dim} \mathcal V = \prod\limits_{i=1}^d N_i$ which depends exponentially on
the dimension $d$.
Then using these basis of the subspace $\mathcal V$ to do the discretization
with the standard process of finite element method will leads to the linear algebraic equation
with the degree of freedom $\prod\limits_{i=1}^d N_i$ which is an exponential function of the dimension $d$.
This is the so called CoD.

In order to overcome the CoD, we use the tensor way to represent the trial function in this paper.
The tensor way
can reduce the complexity of the trail function representation. In the next section, the tensor way also
brings the possibility to do the direct numerical integration for the high dimensional trial functions.
Compared with the methods with Monte-Carlo integration, the tensor
type discretization with the direct numerical integration can improve the stability and accuracy for solving high dimensional partial differential equations.
In this paper, we use the tensor product way to represent the trial functions in the discrete space $\mathcal V$
as follows
\begin{eqnarray}
v(x) = \sum_{j_1=1}^{N_1}\cdots\sum_{j_d=1}^{N_d}T_{j_1,\cdots,j_d}\phi_{i,j_1}(x_1)\cdots\phi_{i,j_d}(x_d),
\end{eqnarray}
where $T = (T_{j_1,\cdots,j_d})_{1\leq j_i\leq N_i} \in \mathcal R^{N_1\times \cdots\times N_d}$
is a tensor with the order of $d$.
There exists the following CP decomposition for the tensor $T$
\begin{eqnarray}\label{CPDecomposition}
T = \sum_{k=1}^P d_k v_{k,1}\otimes \cdots\otimes v_{k,d},
\end{eqnarray}
where $P$ denotes the rank of the function $v(x)$ and $v_{k,i} = (v_{k,i,j})_{1\leq j\leq N_i}$
is a vector of dimensional $N_i$, $1\leq i\leq d$.

Based on the CP decomposition (\ref{CPDecomposition}) of the tensor $T$ and the basis for each
space $\mathcal V_i$, we can use the following way to represent the function $v(x)\in \mathcal V$
\begin{eqnarray}
v(x) = \sum_ {k=1}^P \prod_{i=1}^d \left(\sum_{j=1}^{N_i}v_{k,i,j}\phi_{i,j}(x_i)\right),
\end{eqnarray}
where the coefficient $v_{k,i,j}$ is the $k$-th components in the vector $v_{k,i}$ ($1\leq k\leq P$, $1\leq i \leq d$).

Based on this type of representation, we can build the discrete trial function for
the high dimensional partial differential equation in the following way
\begin{eqnarray}\label{CP_Version}
u_N(x) = \sum_ {k=1}^P \prod_{i=1}^d \left(\sum_{j=1}^{N_i}C_{k,i,j}\phi_{i,j}(x_i)\right),
\end{eqnarray}
where $C_{k,i,j}$ denotes the unknown coefficients for the approximate solution. 
It is worth to mentioning that we do not generate the high dimensional mesh and high dimensional function 
explicitly. We only need to use the one dimensional basis function to build the high dimensional trial functions 
with the tensor product way.

This formula (\ref{CP_Version}) provides a way to reduce the complexity for representing the
trial function for  the high dimensional partial differential equations.  It is easy to know that the
number of coefficients in (\ref{CP_Version}) is
\begin{eqnarray}\label{DOF}
P\left(\sum\limits_{i=1}^dN_i\right).
\end{eqnarray}
Then, if the rank $P$ of the function $u_N(x)$ is polynomial scale of the dimension $d$, the number of
coefficients $C_{k,i,j}$ does not depend exponentially on the dimension $d$. We will also find the decomposition (\ref{CP_Version}) brings the possibility to do the direct numerical integration for these functions, which will be
discussed in the next section.

The aim here is to give a framework to build the trial function for the high dimensional partial differential equations.
With this framework, we can choose different type of basis functions $\phi_{i,j}(x_i)$ to construct the trial function as in (\ref{CP_Version}).
For example, finite element basis, spectral basis, wavelet basis can all be chosen here. Furthermore, the neural network
function is also an important candidate to build the trial function for the high dimensional partial differential equation.
This will be discussed in our another paper.

\section{Solving high-dimensional boundary value problem}\label{Section_BVP}
In this section, we will show the first application of the tensor type discretization method. Here we
are concerned with solving the
high dimensional second order elliptic partial differential equation with homogeneous
Neumann boundary condition: Find $u\in H^1(\Omega)$ such that
\begin{eqnarray}\label{Neumann_Equation}
\left\{
\begin{array}{rcl}
-\Delta u+\pi^2u&=&2\pi^2\sum\limits_{k=1}^d\cos(\pi x_k),\ \ \ x\in \Omega,\\
\frac{\partial u}{ \partial \mathbf n}&=&0,\ \quad\quad\quad\quad\quad\quad\ \ \ x\in\partial\Omega,
\end{array}
\right.
\end{eqnarray}
where the computing domain $\Omega=[0,1]^d$. Then the exact solution is
\begin{eqnarray*}
u(x)=\sum_{k=1}^d\cos(\pi x_k).
\end{eqnarray*}

It is well known that the problem (\ref{Neumann_Equation}) is equivalent to the following  minimization problem
%The solution of this partial differential equation consist with the well-known minimization problem of energy functional as follows
\begin{eqnarray}\label{Optimization_Problem_BV}
u(x)&=&\arg\min_{v\in H^1(\Omega)} I(v)\nonumber\\
&:=&\arg\min_{u\in H^1(\Omega)}\left[\frac{1}{2}\int_{\Omega}\left(|\nabla v|^2+\pi^2v^2\right )dx-\int_{\Omega}2\pi^2\sum_{k=1}^d\cos(\pi x_k)vdx\right].
\end{eqnarray}

In order to solve the optimization  problem (\ref{Optimization_Problem_BV}), we build the tensor type of trial
function $u_N(x)$ in (\ref{CP_Version}),  and define the set of all possible values of coefficients $C_{k,i,j}$.

The trial function set $V$ is modeled by the parameters $C_{k,i,j}$, which take all the possible values
and it is easy to satisfy the condition $V\subset H_0^1(\Omega)$.
The approximate solution  $u_N(x)$ for (\ref{Optimization_Problem_BV}) can be defined as the 
solution of the following optimization problem
\begin{eqnarray}\label{approx_opt}
u_N(x)&:=&\arg\min_{u_N\in V}\left[\frac{1}{2}\int_{\Omega}\left(|\nabla v|^2+\pi^2v^2\right )dx-\int_{\Omega}2\pi^2\sum_{k=1}^d\cos(\pi x_k)vdx\right].
\end{eqnarray}
Note that the terms of (\ref{approx_opt}) includes the high dimensional integration.
It is well known that the high dimensional integration has the CoD.
Always, the Monte-Carlo method is used to do these high dimensional integration. Fortunately, since the tensor
decomposition structure of the trail function in this paper, we can do the direct high dimensional integration without Monte-Carlo process.

From the definition of $u_N(x)$ in (\ref{CP_Version}), we have the following representation for $u_N^2(x)$
\begin{eqnarray}\label{u_N_2}
u_N^2(x) &=&\sum_{k=1}^P \prod_{i=1}^d\left(\sum_{j=1}^{N_i} C_{k, i, j} \phi_{i, j}\left(x_i\right)\right)
\sum_{\ell=1}^P \prod_{m=1}^d\left(\sum_{n=1}^{N_m} C_{\ell, m, n} \phi_{m, n}\left(x_m\right)\right) \nonumber\\
&=&\sum_{k=1}^P \sum_{\ell=1}^P \prod_{i=1}^d\left(\sum_{j=1}^{N_i} C_{k, i, j} \phi_{i, j}\left(x_i\right)\right)\left(\sum_{n=1}^{N_i} C_{\ell, i, n} \phi_{i, n}\left(x_i\right)\right) .
\end{eqnarray}
Then the following quadrature scheme holds
\begin{eqnarray}\label{Quadrature_1}
\int_\Omega u_N^2(x)dx &=&\sum_{k=1}^P \sum_{\ell=1}^P \prod_{i=1}^d\int_{\Omega_i}
\left(\sum_{j=1}^{N_i} C_{k, i, j} \phi_{i, j}\left(x_i\right)\right)
\left(\sum_{n=1}^{N_i} C_{\ell, i, n} \phi_{i, n}\left(x_i\right)\right) dx_i\nonumber\\
&=&\sum_{k=1}^P \sum_{\ell=1}^P \prod_{i=1}^d
\left(\sum_{j=1}^{N_i} \sum_{n=1}^{N_i}  C_{k, i, j}
C_{\ell, i, n}  \int_{\Omega_i}\phi_{i, j}\left(x_i\right) \phi_{i, n}\left(x_i\right)dx_i\right).
\end{eqnarray}
From (\ref{Quadrature_1}), in order to do the high dimensional integration $\int_\Omega u_N^2(x)dx$, we only need
to compute the one dimensional integration $ \int_{\Omega_i}\phi_{i, j}\left(x_i\right) \phi_{i, n}\left(x_i\right)dx_i$
for $1\leq j,n\leq N_i$ and $1\leq i\leq d$.
Other terms in (\ref{approx_opt}) can be treated with the similar way in (\ref{Quadrature_1}).

\begin{theorem}\label{Theorem_Gauss}
Assume that the function $u(x)$ is defined as (\ref{CP_Version}) in the $d$-dimensional tensor product domain $\Omega$.
Let $N:=\max\{N_1, \cdots, N_d\}$ and $T_1$ denote the computational complexity
for the following $1$-dimensional integration
\begin{eqnarray}
\int_{\Omega_i}\phi_{i,j}(x_i)\phi_{i,n}(x_i)dx_i,\ \ \  1\leq j,n\leq N_i,\ i=1, \cdots, d.
\end{eqnarray}
The computational complexity for the numerical integration (\ref{Quadrature_1}) can be bounded by
$\mathcal O\big(dT_1P^2N\big)$, which is the polynomial scale of the dimension $d$.
\end{theorem}
%\begin{proof}
%First, we show that the number of $j_{\beta,\ell}$ in the last summation of (\ref{eq_decomposition_prod})
%is no more than $k$.
%This result can be easily proved by the following inequality
%\begin{eqnarray}
%\sum_{\beta\in\mathcal B_\alpha}\alpha_\beta=|\alpha|\leq k.
%\end{eqnarray}
%Then, by direct calculation, the computational complexity of (\ref{eq_I_tensor_form})
%can be bounded by $\mathcal O\big(dT_1 p^2 \big)$.
%This is the desired result and the proof is complete.
%\end{proof}

\begin{remark}
For some special equations, we can choose the basis with orthogonal property such that the computational 
complexity for the integration in Theorem \ref{Theorem_Gauss} can be reduced to $\mathcal O\big(dT_1PN\big)$. 
\end{remark}

If $u_N(x)$ does not have the tensor form, when using the same quadrature scheme, it is easy to know that
the computational complexity for the high dimensional integration
depends exponentially on the dimension.
% defined by the NN is $\mathcal O\big((dqT_1+kT_d)((d+m)^m+k)^kN^d\big)$,
%where $T_d$ denotes the complexity for the $d$-dimensional function evaluation operations.

In this paper, the gradient descent (GD) method is adopted to solve the optimization problem (\ref{approx_opt}).
The GD scheme can be described as follows:
\begin{eqnarray}\label{gd}
\theta^{(k+1)}=\theta^{(k)}-\eta\nabla \mathcal R(u^{(k)}(x)),
\end{eqnarray}
where $\theta^{(k)}$ denotes the parameters after the $k$-th GD step, $\eta$ is the learning rate (step size). 
Actually, in this paper, we use the alternating coordinate gradient descent method to solve the optimization problem (\ref{approx_opt})
since in each dimensional, the target function is an quadratic form and it is can be easily optimized.

%It is well known that the general FNN-based machine learning method always uses
%Monte-Carlo procedure  to do the high dimensional integration. With the tensor structure, we can do the direct numerical integrations for the concerned trial functions in this paper.
%For example, we can use the fixed quadrature points (i.e. Gauss points)
%to do the numerical integration.
%If we use the fixed quadrature points for FNN functions, the computational work for the numerical integration will depend exponentially on the dimension $d$.
%In order to avoid the  CoD, this is the reason to use the Monte-Carlo method and
%the stochastic gradient descent (SGD) method \cite{KingmaAdam}  in the numerical implementation for solving the high-dimensional PDEs by FNN-based method \cite{EYu}.

Fortunately, based on the tensor structure in this paper, Theorem \ref{Theorem_Gauss} shows that
the high dimensional numerical integration does not encounter CoD %``curse of dimensionality''
since the computational work can be bounded by the polynomial scale of dimension $d$.
Due to the high accuracy of high dimensional integration and then the target function,
the high accuracy of the tensor type discretization method can be guaranteed.

\section{Numerical experiments}\label{Section_Numerical}
In this section, we come to investigate the numerical performance of the tensor type discretization method for the 
 problem  (\ref{Neumann_Equation}).   The aim here is to check the efficiency and accuracy of the proposed 
 numerical method in this paper.

%In each example we will check the performances of tensor type discretization for three cases of dimensions:
%$d=128$, $d=256$ and $d=512$. The solution $u(x)$ is approximated by the tensor structure
%$u(x)$ which is defined by (\ref{CP_Version}).
For the simplicity, the linear finite element basis functions \cite{Braess} 
are adopted to act as  the basis $\phi_{i,j}(x_i)$ for
$1\leq j\leq N_i$ and $1\leq i\leq d$ for the tensor type of trial functions  (\ref{CP_Version}).  
Since the exact solution $u$ is known, we can check the exact error for the approximate solution with the help 
of the numerical integration for the tensor type functions.  
In order to show the convergence behavior and accuracy of approximations by tensor type
discretization, we define the $H^1(\Omega)$ and  $L^2(\Omega)$ error estimates as follows 
\begin{eqnarray*}
{\rm Err}_{H^1} := \frac{\|u(x)-u_N(x)\|_1}{\|f\|_0},\ \ \ {\rm Err}_{L^2} := \frac{\|u(x)-u_N\|_0}{\|f\|_0}.
\end{eqnarray*}
In all tests, we set the rank $P=2d$. All the numerical tests are done by the laptop with 
I7 CPU (2.8G) and 8G memory and the code is written with Python.  
The optimization problem (\ref{approx_opt}) is solved with the alternating coordinate gradient descent method.

\begin{figure}[htb]
\centering
\includegraphics[width=7cm,height=5cm]{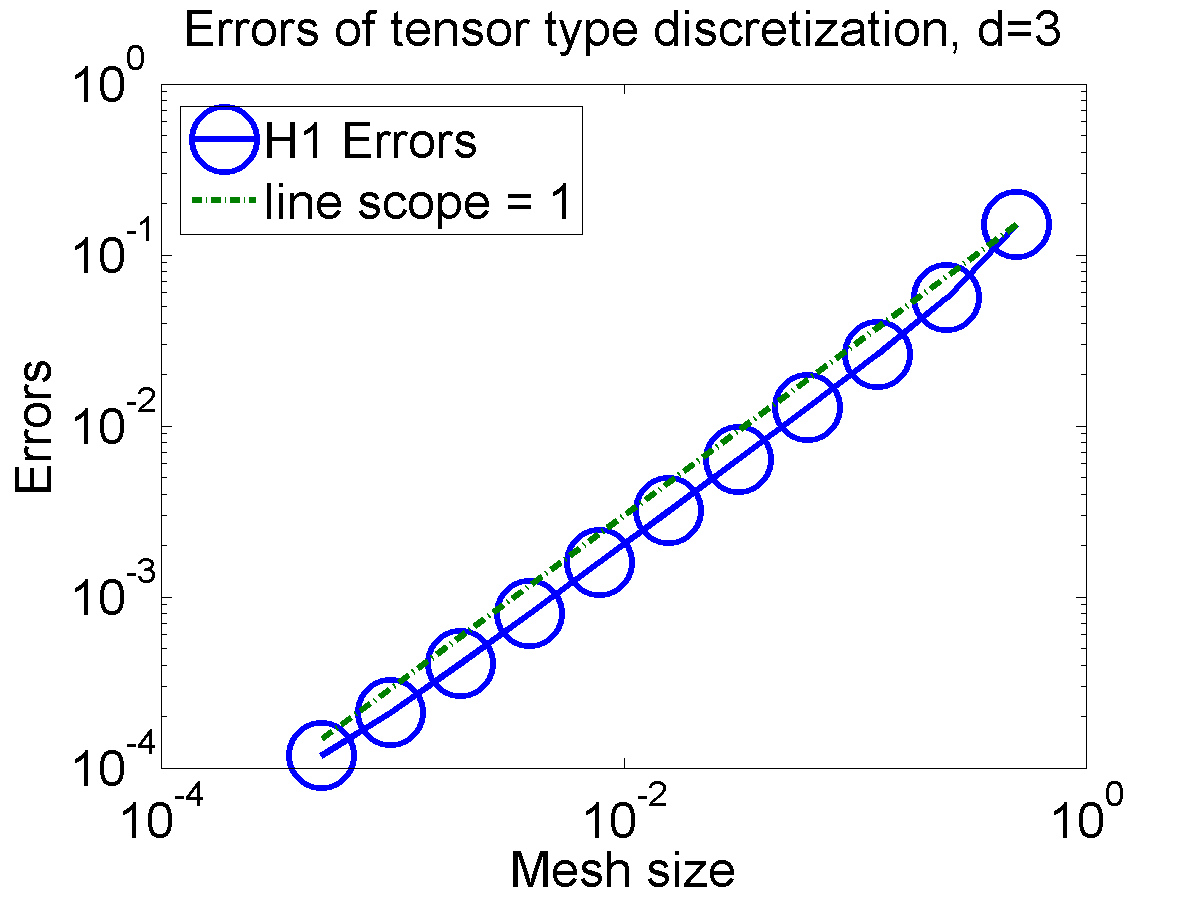}
\includegraphics[width=7cm,height=5cm]{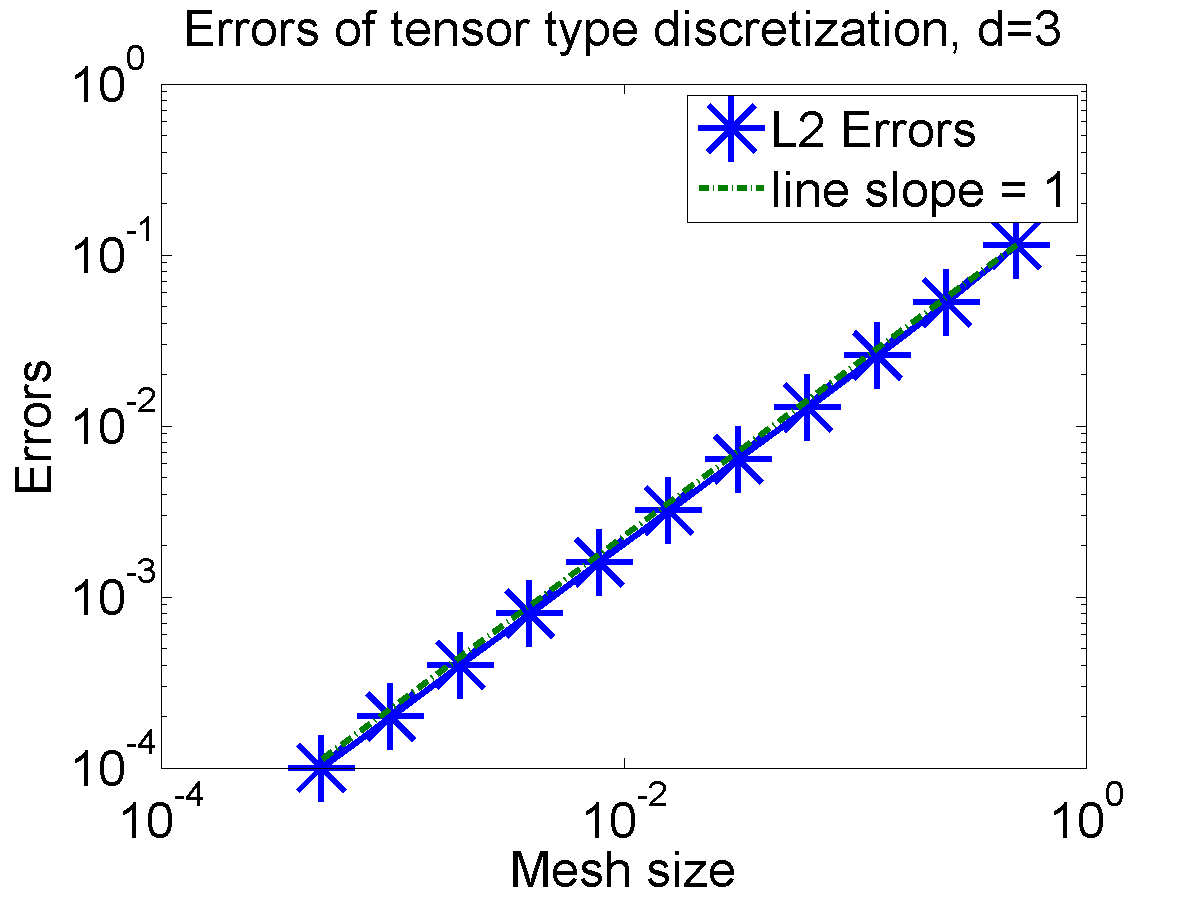}
\caption{Errors for the three dimensional case with the tensor type of discretization method 
with the sequence of mesh sizes $h = [1/2, 1/4, 1/8, 1/16, 1/32,1/64,1/128, 1/256,1/512,1/1024, 1/2048]$.}\label{fig_3d}
\end{figure}

In order to show the comparison with the standard finite element method, we fist show the numerical results for the 
three dimensional case $d=3$. The corresponding numerical results are shown in Figure \ref{fig_3d}. 

From Figure \ref{fig_3d}, we can find that the tensor type of discretization can also get the optimal convergence 
order for the $H^1(\Omega)$ norm errors. Different form the standard finite element method, the $L^2(\Omega)$ errors 
can also only have the first convergence order.  As we know,  in this case, we can also use a finite element packages 
 to solve the three dimensional boundary value problem (\ref{Neumann_Equation}). But in order to obtain the same accuracy, 
 we need to generate the three dimensional mesh with the mesh size $h=1/2048$. Then the corresponding number of 
 tetrahedral is $6\times 2048^3 \approx 5.15\times 10^{10}$ and we must use the high performance computer and the parallel 
 method.  For the comparison, with the tensor type of discretization, we only need to use a laptop to solve 
 the same problem. 
 
\begin{figure}[htb]
\centering
\includegraphics[width=7cm,height=5cm]{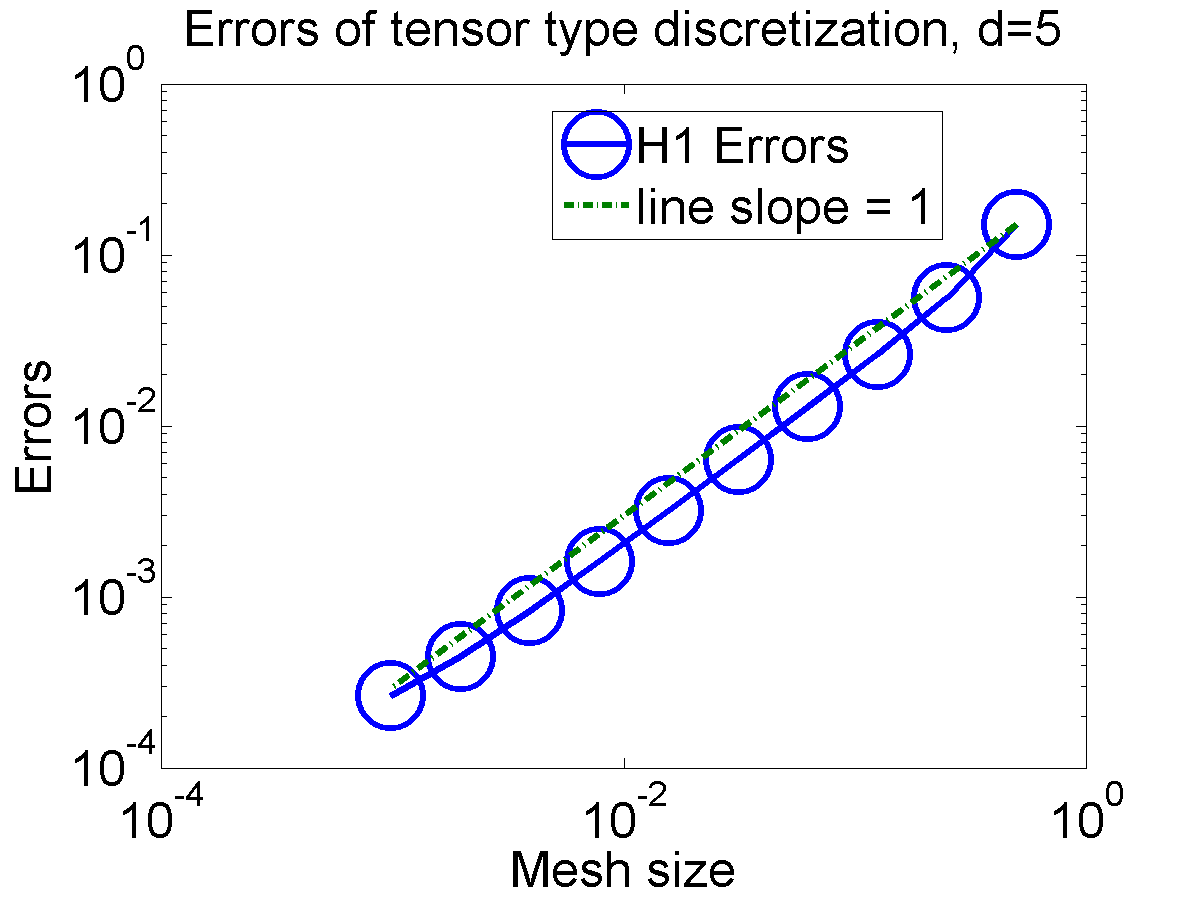}
\includegraphics[width=7cm,height=5cm]{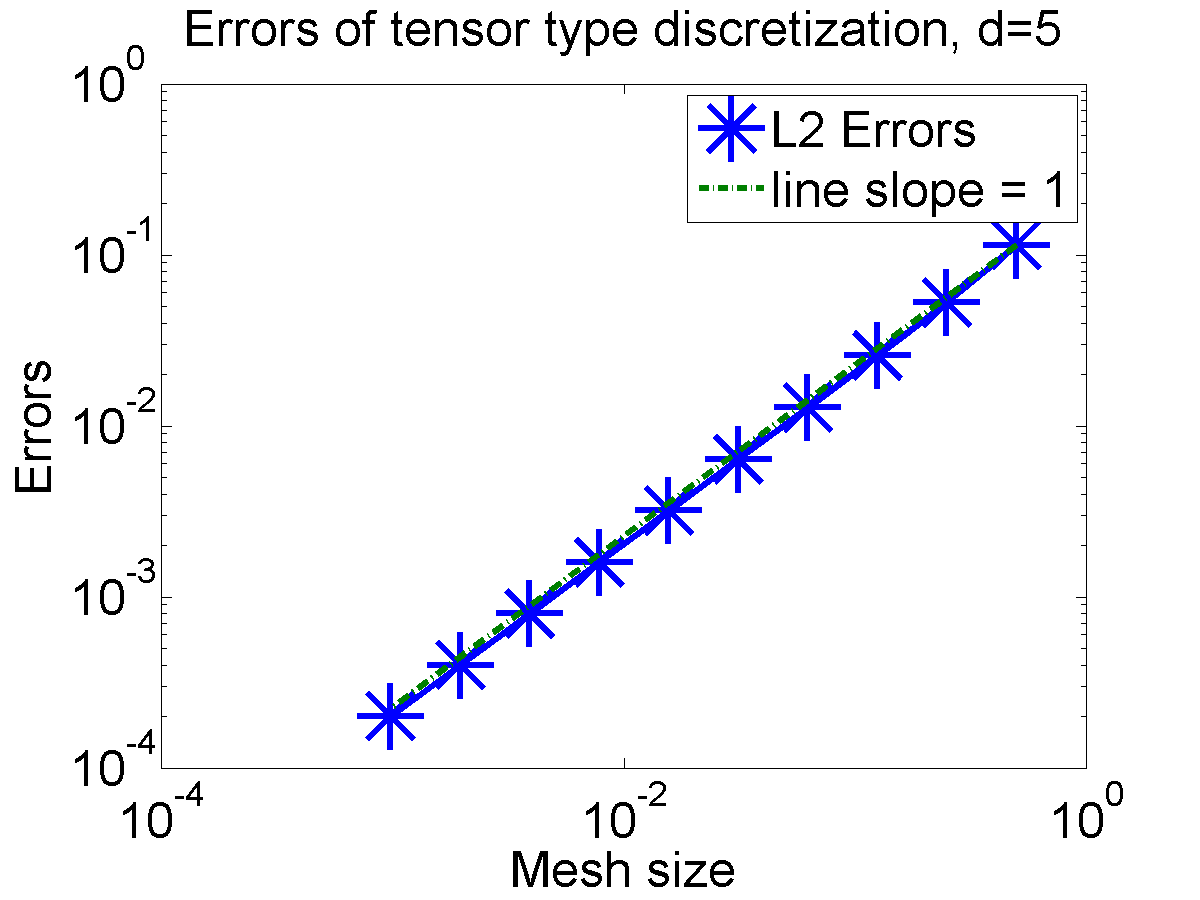}
\caption{Errors for the five dimensional case with the tensor type of discretization method
with the sequence of mesh sizes $h = [1/2, 1/4, 1/8, 1/16, 1/32,1/64,1/128, 1/256,1/512,1/1024]$.}\label{fig_5d}
\end{figure}

 In this section, we also do the numerical experiments for the five dimensional boundary value problem 
 (\ref{Neumann_Equation}).  The corresponding numerical results are shown in Figure \ref{fig_5d}. 

From Figure \ref{fig_5d}, we can also find that the tensor type discretization method can obtain the 
optimal convergence order for the $H^1(\Omega)$ error estimates and the $L^2(\Omega)$ errors also has only 
first order convergence.

\section{Conclusions}
In this paper, we present the tensor type of discretization method for solving high dimensional partial 
differentia equations.  
Since the proposed discretization method has the tensor product structure, we can do the direct numerical integration
without the help of Monte-Carlo process for the high dimensional trial functions and their inner-products.
The direct numerical integration for high dimensional trial functions can improve the stability and accuracy for
solving high dimensional partial differential equations.
We believe the ability of direct numerical integration will bring more applications.

Based on the constructing way and numerical integration method in this paper, we can do the following extensions:
\begin{enumerate}
\item Different types of basis functions can be adopted to build the trial functions for solving
 high dimensional partial differential equations.

\item TThe neural network can be used to construct the tensor type discretization here. 
This means the neural network can be combined with the tensor type discretization to build a new type of
machine learning method for solving high dimensional partial differential equations.

\item Since the direct numerical integration for the trial functions, we can do the posteriori error estimate 
for the approximate solutions \cite{Braess}.  This also provide an adaptive way to chose the rank for building the trial functions.

\item Since we need to solve  the nonlinear optimization problem for the tensor type  discretization 
method, the augmented subspace method \cite{Xie,Xie2} can be used to improve the efficiency of solving the 
nonlinear optimization problems. 

\end{enumerate}
More applications to other types of problems will also be considered in the future.

\end{document}